\documentclass[12pt,reqno]{article}

\usepackage[usenames]{color}
\usepackage{amssymb}
\usepackage{amsmath}
\usepackage{amsthm}
\usepackage{amsfonts}
\usepackage{amscd}
\usepackage{graphicx}

\DeclareMathOperator{\Li}{Li}

\usepackage[colorlinks=true,
linkcolor=webgreen,
filecolor=webbrown,
citecolor=webgreen]{hyperref}

\definecolor{webgreen}{rgb}{0,.5,0}
\definecolor{webbrown}{rgb}{.6,0,0}

\usepackage{color}
\usepackage{fullpage}
\usepackage{float}

\usepackage{graphics}
\usepackage{latexsym}
\usepackage{epsf}
\usepackage{breakurl}

\def\gcd{\mathrm{gcd}}
\def\Li{\mathrm{Li}}

\def\modd#1 #2{#1\ \mbox{\rm (mod}\ #2\mbox{\rm )}}

\begin{document}

\begin{center}
\epsfxsize=4in

\end{center}

\theoremstyle{plain}
\newtheorem{theorem}{Theorem}
\newtheorem{corollary}[theorem]{Corollary}
\newtheorem{lemma}[theorem]{Lemma}
\newtheorem{proposition}[theorem]{Proposition}

\theoremstyle{definition}
\newtheorem{definition}[theorem]{Definition}
\newtheorem{example}[theorem]{Example}
\newtheorem{conjecture}[theorem]{Conjecture}

\theoremstyle{remark}
\newtheorem{remark}[theorem]{Remark}

\begin{center}
\vskip 1cm{\LARGE\bf  Representing positive integers as a sum of a squarefree number and a small prime
}

\vskip 1cm
\large
Jack Dalton and Ognian Trifonov\\
University of South Carolina\\
Department of Mathematics \\
Columbia, SC 29208 \\
\href{mailto:jrdalton@math.sc.edu}{\tt jrdalton@math.sc.edu} \\
\href{mailto:trifonov@math.sc.edu}{\tt trifonov@math.sc.edu} \\
\end{center}

\vskip .2 in

\begin{abstract}
We prove that every positive integer $n$ which is not  equal to $1$, $2$, $3$, $6$, $11$, $30$, $155$,  or $247$ can be represented as a sum of a squarefree number and a prime not exceeding $\sqrt{n}$.
\end{abstract}

\section{Introduction}

A positive integer is squarefree if it is not divisible by the square of a prime number. In 1931 Estermann \cite{Est} obtained an asymptotic formula for the number of representations of a positive integer as the sum of a squarefree number 
and a prime. As a consequence, he showed that every sufficiently large integer is the sum of a squarefree number 
and a prime.

In 2017 Dudek \cite{Dud} showed that every integer greater than two is a sum of a squarefree number 
and a prime. 
 
In 1935,  Erd\H{o}s \cite{Erd} showed that every sufficiently large integer $ \not \equiv \modd{1} {4}$ can be represented as a sum of a square of a prime and a squarefree integer. The condition $n \not \equiv \modd{1} {4}$ is needed since if $p$ is an odd prime, then $p^2 \equiv \modd{1} {4}$. Therefore,  if $p$ is an odd prime and $n  \equiv \modd{1} {4}$, then $4 \mid n-p^2$, and $n-p^2$ is not squarefree. The result of Erd\H{o}s was completed by Dudek and Platt \cite{DudPl} who proved in 2016 that every integer $n \geq 10$ such that $n \not \equiv \modd{1} {4}$ can be represented as the sum of a square of a prime and a squarefree integer.

Thus, one can get a complete answer to the question of for which positive integers $n$ there exists a prime $p < \sqrt{n}$ such that $n-p^2$ is squarefree.

A related question is for which positive integers $n$, there exists a prime $p < \sqrt{n}$ such that $n-p$ is squarefree.
We answer the latter question completely.
 
\begin{theorem} \label{sqsp}
Every positive integer $n$ can be represented in the form $n=s+p$ where $s$ is a squarefree number and $p$ is a prime with $p \leq \sqrt{n}$, except when $n \in \{1, 2, 3, 6, 11, 30, 155, 247\}$.
\end{theorem}

 Since $2^2 \mid 6-2$; $3^2 \mid 11-2$, $2^3 \mid 11-3$; $2^2 \mid 30-2$, $3^3 \mid 30-3 $, $5^2 \mid 30-5 $; $3^2 \mid 155-2$,  $2^3 \mid 155-3$, $5^2 \mid 155-5$, $2^2 \mid 155-7$, $2^2 \mid 155-11$; and $7^2 \mid 247-2$, 
 $2^2 \mid 247-3$, $11^2 \mid 247-5$, $2^2 \mid 247-7$, $2^2 \mid 247-11$, and $3^2 \mid 247-13$ we see that indeed the integers $n   \in \{1, 2, 3, 6, 11, 30, 155, 247\}$ cannot be represented as a sum of a squarefree number and a prime not exceeding $\sqrt{n}$. 
 
A natural question is can $ \sqrt{n}$  in Theorem \ref{sqsp} be replaced by $n^{\theta}$ for $\theta < 1/2$. 
 
In 2015, Filaseta, Graham, and the second author wrote a paper \cite{FGT} which considered distribution of various arithmetic sequences. In particular, replacing $p$ by $-p$ in Theorem~5.1 of \cite{FGT} one obtains the following theorem.
 
 \begin{theorem}[Filaseta, Graham, and T.] \label{FGT thm}
 
 There exist effectively computable constants $C_0$ and $n_0$ such that for each integer $n \geq n_0$ at least one-fifth of the primes $p \leq C_0 n^{1/5} \log ^2 n$ are such that $n-p$ is squarefree.
 
 \end{theorem}
 
 The authors of the paper \cite{FGT} did not compute the  values of $C_0$ and $n_0$. However,  it is certain that it will be impossible to check which integers $n \leq n_0$ do not satisfy Theorem \ref{FGT thm} by direct computation.
 
 An easy corollary of the above theorem is that for each $\theta > 1/5$ there exists an effectively computable constant $n_{\theta }$ such that each integer $n \geq n_{\theta }$ can be represented as a sum of a squarefree number and a prime not exceeding $n^{\theta}$. At present it appears that to get below the exponent $1/5$ one will need to improve the gap result about squarefree numbers.
 
 The proof of our main result follows the method of proof of Theorem \ref{FGT thm}. However, to obtain the result for all but eight positive integers we substantially sharpened the estimates in each step and supplemented the proof with nontrivial amount of computations. 
 
Next, we outline the proof of Theorem \ref{sqsp}.  First, we verify the theorem for all $n \leq 10^9$. This computation allows us to obtain the first few terms of the following sequence. 

Let $b_k$ be the least positive integer $n > p_l$ such that none of the integers $n - p_1, \ldots, n - p_k$ are squarefree. Here and throughout the paper $p_l$ denotes the $l$-th prime number. 

We have $b_1=6$, $b_2 = 11$, $b_3 = 30$, $b_4=b_5=155$, $b_6=247$, $b_7=5753$, $b_8=b_9=b_{10}=b_{11}=90263$, $b_{12}=1481287$, $b_{13}=b_{14}=b_{15}=7409327$. Also, $b_{16}> 10^9$.

Clearly, ${\displaystyle b_n < \prod_{l=1}^n p_l^2 }$ for $n>1$, since by the Chinese Remainder Theorem, the system of congruences $n \equiv p_l \pmod{p_l^2}$ has a positive solution not exceeding ${\displaystyle \prod_{l=1}^n p_l^2 }$. 

One can do much better. For example, if we want to find $n > 53=p_{16}$ such that $n - p_l$ is not squarefree for $l=1,\ldots,16$, one can pick $n \equiv 3 \pmod{4}$. Then $4 | n - p$ for $p \in \{3,7,11,19,23,31,43,47 \}$. Also, 
if $n \equiv 8 \pmod{9}$, then  $9 | n - p$ for $p \in \{17,53 \}$. After that, we pick $n \equiv 2 \pmod{5^2}$, $n \equiv 5 \pmod{7^2}$, $n \equiv 13 \pmod{11^2}$, $n \equiv 29\pmod{13^2}$, $n \equiv 37 \pmod {17^2}$, and
$n \equiv 41 \pmod{19^2}$. Solving the last system of congruences, we obtain that its least positive solution is  $23708451225527$, thus $b_{16} \leq 23708451225527$. 

From now on we assume $n > 10^9$. We show that for such $n$, there exists a prime $p \leq  \sqrt{n}$, such that $n-p$ is squarefree. 

There are $\pi(\sqrt{n})$ primes which do not exceed $\sqrt{n}$. We will show that there are less than $\pi(\sqrt{n})$ primes $p \leq \sqrt{n}$ such that $n-p$ is not squarefree. Establishing the last statement,  will imply that for some prime $p \leq \sqrt{n}$, 
$n - p = s$ where $s$ is squarefree completing the proof of the theorem. 

If for some prime $p \leq \sqrt{n}$, $n-p$ is not squarefree, then $q^2 | n-p$ for some prime $q$, that is 
\begin{equation}\label{cong}
  p \equiv \modd{n} {q^2}.
\end{equation}

Note that congruence \eqref{cong} can hold only for $q < \sqrt{n}$, since $0 < n-p < n$. Furthermore, 
if $q | n$ the congruence \eqref{cong} has at most one solution, when $p=q$. Also, then the number of solutions 
of \eqref{cong} is $\pi(\sqrt{n};q^2,n)$. 

Above we have used the standard notation $\pi(a;b,c)$ to denote the number of primes not exceeding $a$ which are in the arithmetic sequence $x \equiv \modd{c} {b}$. 

We estimate the number of primes $p$ such that congruence \eqref{cong} holds for some prime $q$ in different ways depending on the size of $q$. For $q \in \{ 2,3,5,7 \}$ we use the paper  \cite{BenMar} 
of Bennett et. al. on explicit bounds for primes in arithmetic sequences. For $q \in [11, n^{1/8}]$, we use a version of the Brun-Titchmarsh Theorem due to Montgomery and Vaughan \cite{MV}. For $q \in [n^{1/8}, \sqrt{n}/(c\log n)]$ we use elementary bounds, and for $n \in  [\sqrt{n}/(c\log n), \sqrt{n})$ we use an estimate based on first differences of values of a function. Combining the above estimates we were able to show that the theorem holds for $n > 59^8$. Next, considering 
$n \in [p_k^8, p_{k+1}^8]$ for $k=6,\ldots, 16$ we were able to get even sharper bounds for certain sums over primes which appear in our estimates and establish the theorem in the last $11$ intervals. Since $p_6=13$ and 
$13^8 < 10^9$ this completes the proof.

\section{Proof of the main result}

\begin{proof}
First, by direct computation we confirmed the theorem for $n \leq 10^9$. To save memory we check the result in intervals of the form $I_k :=[k\cdot 10^7, (k+1)\cdot 10^7]$. Computing the indicator function of squarefree integers in $I_k$ is relatively fast since for each prime $p$ the number of multiples of $p^2$  in $I_k$ is at most $10^7/p^2 + 1$ and one only needs to check multiples of $p^2$ with $p \leq \sqrt{(k+1)\cdot 10^7}$. Checking the theorem in $I_k$ is also fast because for most $n$, $n-2$
is squarefree and for each $n \leq 10^9$ it took no more than $17$ checks to find a prime $p$ such that $n-p$ is squarefree. Using python on a laptop with an AMD Ryzen 5 3500U processor and 8MB of memory, checking the theorem on each $I_k$ with 
$k \in [1,99]$ took less than a minute, so the whole computation was complete in under two hours.
 
From now on, we assume $n > 10^9$.

 We will use the notation $ \sideset{}{'}  \sum $ to indicate a sum over prime numbers only. For example, ${\displaystyle   \sideset{}{'} \sum_{q \leq x} \  1}$ is the number of prime numbers up to $x$, that is $\pi (x)$.

 We  estimate the number of primes $p < \sqrt{n}$ satisfying congruence \eqref{cong} in different ways depending on the size of $q$.

{\bf Case 1.} $q \in \{2,3,5,7 \}$.

Let $Q_1$ be the number of primes less that $\sqrt{n}$ which satisfy congruence \eqref{cong} for some $q \in \{2,3,5,7 \}$.

Here we use the following theorem which is one of the statements in Corollary 1.7 from the paper \cite{BenMar}.

\begin{theorem} (M.~Bennett, G.~Martin, K.~O'Bryant, and A.~Rechnitzer) \label{BM}

Let $a$ and $q$ be integers with $1 \leq q \leq 10^5$ and $\gcd (a,q)=1$. If $x \geq 10^6$, then
$$ \left | \pi(x;q,a) - \frac{\Li (x)}{\varphi(q)} \right | < 0.027 \frac{x}{\log ^2 x}.$$
\end{theorem}
   
 Above, by $\Li (x)$ the authors  mean the function defined by ${\displaystyle \Li (x) = \int_2^x \frac{dt}{\log t}}$. Also, everywhere in this paper $\log x$ means $\log_e x = \ln x$.

  To be able to use the above theorem, we assume $n \geq 10^{12}$ and we will deal with the case $10^9 < n < 10^{12}$ later. 
  
  Using inclusion-exclusion for $q=2$ and $q=3$, and that congruence \eqref{cong} has $    \pi(\sqrt{n};q^2,n)$ solutions for fixed $q$, we obtain, 
  
  \begin{equation} \label{Q11}
  Q_1 \leq     \pi(\sqrt{n};4,n)+    \pi(\sqrt{n};9,n) -  \pi(\sqrt{n};36,n) + \pi(\sqrt{n};25,n) + \pi(\sqrt{n};49,n).
  \end{equation}
  
  When $\gcd(n,210)=1$, Theorem \ref{BM} applies, and we obtain
  
  \begin{equation} \label{Q12}
  Q_1  \leq \Li (\sqrt{n})\left ( \frac{1}{2} + \frac{1}{6} - \frac{1}{12} + \frac{1}{20} + \frac{1}{42} \right ) + 0.54\frac{\sqrt{n}}{\log ^2 n}.
  \end{equation}

 The following is Lemma 5.9 of the paper \cite{BenMar}.

\begin{lemma} \label{Li(x)}
For  $ x \geq 1865$, 
\begin{equation} \label{Li}
\Li (x) < \frac{x}{\log x} \left (  1 + \frac{3}{2\log x}  \right )  \quad . 
\end{equation}
\end{lemma}
   
 Using the above lemma and \eqref{Q12} we get 
 
 \begin{equation} \label{Q13}
  Q_1  < \frac{46\sqrt{n}}{35\log n}\left (  1 + \frac{3}{\log n}    \right ) + 
0.54\frac{\sqrt{n}}{\log ^2 n}\ = \  \frac{46\sqrt{n}}{35\log n} +  \frac{1569\sqrt{n}}{350\log^2 n}.
\end{equation}   
 
If $7 | n$, $\pi(\sqrt{n};49,n) \leq 1$ and we need to replace to upper bound for $\pi(\sqrt{n};49,n) \leq 1$ used in \eqref{Q13} which is 
${\displaystyle  \frac{ \sqrt{n}}{21\log n}\left (  1 + \frac{3}{\log n}     \right ) + 
0.108\frac{\sqrt{n}}{\log ^2 n}}\ $ by $1$.

However, the function ${\displaystyle f(n) = \frac{\sqrt{n}}{\log n} }$ is increasing for $n > e^2$ and $f(10^{12})=36192.2\ldots$, so 
$$   \frac{ \sqrt{n}}{21\log n}\left (  1 + \frac{3}{\log n}      \right ) + 
0.108\frac{\sqrt{n}}{\log ^2 n} > 36192/21 > 1,$$
for $n \geq 10^{12}$, so equation \eqref{Q13} holds when  $\gcd(n,30)=1$ (regardless of whether $7 \mid n$ or $7  \nmid n$). 

We argue similarly when $5 \mid n$ to conclude that equation \eqref{Q13} holds when $\gcd(n,6)=1$ and $n \geq 10^{12}$. 

When $3 \mid n$ we need to replace to upper bound for $\pi(\sqrt{n};9,n) - \pi(\sqrt{n};36,n) \leq 1$ used in \eqref{Q13} which is 
${\displaystyle  \frac{ \sqrt{n}}{12\log n}\left (  1 + \frac{3}{\log n}     \right ) + 
0.216\frac{\sqrt{n}}{\log ^2 n}}\ $ by $1$, which is not a problem for $n \geq 10^{12}$. Thus, equation \eqref{Q13} holds 
when $\gcd(n,2)=1$ and $n \geq 10^{12}$. 

We argue similarly when $2 | n$ to show that equation \eqref{Q13} holds for all positive integers $n \geq 10^{12}$.

\vspace{0.25 in}
\noindent
{\bf Case 2.} $11 \leq q \leq n^{1/8}$

Since $n \geq 10^{12}$ we have $n^{1/8} > 11$. 

Let $Q_2$ be the number of primes $p \leq  \sqrt{n}$  such that congruence \eqref{cong} holds for some $q \in [11, n^{1/8}]$. 
Thus, 

\begin{equation} \label{Q21}
Q_2  \leq  \sideset{}{'} \sum_{11 \leq q \leq n^{1/8}}  \  \pi(\sqrt{n};q^2,n).
\end{equation}

In Case 2 the  following theorem from the paper \cite{MV} which  is a version of Brun-Titchmarsh's inequality will be helpful. 

\begin{theorem} (H.~L.~Montgomery and R.~C.~Vaughan) \label{MV}
Let $m > 0$  and $l$  be integers with $(m, l) = 1$, and let $x > m$ be a real number. Then
$$\pi (x ; m , l) \leq \frac{2x}{ \varphi  (m)  \log(x/m) }.$$
\end{theorem}

We apply  Theorem \ref{MV} with $x = \sqrt{n}$, $m=q^2$, and $l = n$. When $q \nmid n$ and 
$q^2 < \sqrt{n}$ the conditions of the theorem hold.  So, suppose that $q \leq n^{1/8}$ and  $q \nmid n$. 
Then,
\begin{equation} \label{c11}
\pi(\sqrt{n};q^2,n) \leq \frac{2\sqrt{n}}{ (q-1)q  \log(\sqrt{n}/q^2) }.
\end{equation}

When $q \mid n$, then $\pi(\sqrt{n};q^2,n) \leq 1$. However, for $q \geq 2$, $ \log(\sqrt{n}/q^2) \leq (\log n)/2$. Also, $q(q-1) < n^{1/4}$. So, 
$$ \frac{2\sqrt{n}}{ (q-1)q  \log(\sqrt{n}/q^2) } > \frac{2n^{1/4}}{\log n}.$$
The function   ${\displaystyle f_1(n) = \frac{n^{1/4}}{\log n} }$ is increasing for $n > e^4$ and $f_1(10^{9})=8.58\ldots$, so equation \eqref{c11} holds when $q \mid n$, as well.

Note that $$\frac{1}{\log (\sqrt{n}/q^2)} = \frac{2}{\log n - 4\log q} = \frac{2}{\log n}
\left ( \frac{1}{1 - \frac{4\log q}{\log n}} \right ).$$

Next, use that $$\frac{1}{1-r} \leq 1+2r$$ for $r \leq 1/2$ with ${\displaystyle r = \frac{4\log q}{\log n}}$ to obtain, 
$$\frac{1}{\log (\sqrt{n}/q^2)} \leq  \frac{2}{\log n} + \frac{16\log q}{\log ^2 n}. $$

Combining the above inequality with  equation \eqref{c11} we get 
$$ \pi(\sqrt{n}; q^2, n) \leq \frac{4\sqrt{n}}{(q-1)q \log n} + \frac{32\sqrt n \log q}{(q-1)q \log ^2 n}.$$

Adding the above inequality for primes $q$ in $[11, n^{1/8}]$ we obtain

\begin{equation} \label{r2}
 \sideset{}{'} \sum_{ 11 \leq q \leq n^{1/8}}   \pi(\sqrt{n}; q^2, n)  \leq c_1(11,n^{1/8})  \frac{4\sqrt{n}}{\log n}+ c_2(11,n^{1/8})  \frac{32\sqrt{n}}{\log^2 n},
\end{equation}
where we define $$c_1(A,B) = \sideset{}{'}  \sum_{  A \leq q \leq  B}   \frac{1}{q(q-1)}, \quad c_2(A,B) =\sideset{}{'}  \sum_{  A \leq q \leq  B}  \frac{\log q}{q(q-1)}.$$  
.

Note that $$c_1(A,B) < C_1(A):=  \sideset{}{'} \sum_{ A \leq q}  \frac{1}{q(q-1)},$$
and
$$c_2(A,B) < C_2(A):=  \sideset{}{'} \sum_{ A \leq q}  \frac{\log q}{q(q-1)}.$$

So, in the case when we have no upper bound on $n$ we will use the estimate
\begin{equation} \label{r21}
 Q_2 \leq \sideset{}{'} \sum_{ 11 \leq q \leq n^{1/8}}  \pi(\sqrt{n}; q^2, n)  \leq    \frac{4C_1(11)\sqrt{n}}{\log n} +    \frac{32C_2(11)\sqrt{n}}{\log^2 n}.
\end{equation}

{\bf Case 3.} $n^{1/8} < q \leq  n^{1/4}$

Let $Q_3$ be the number of primes $p \leq  \sqrt{n}$  such that congruence \eqref{cong} holds for some $q \in (n^{1/8}, n^{1/4}]$. 
Thus, 

\begin{equation} \label{Q31}
Q_3  \leq  \sideset{}{'} \sum_{n^{1/8} < q \leq n^{1/4}}  \  \pi(\sqrt{n};q^2,n).
\end{equation}

Here  we use the following lemma.

\begin{lemma} \label{elsie}  
Let $x > 1$, let $a$ be an integer, and let $q>1$ be a prime $p$, or a power of a prime, $p^k$ with $p>3$. Then,
\begin{equation} \label{721}
  \pi(x;q,a) \leq \frac{x}{3q}+2.
\end{equation}
\end{lemma}

\begin{proof}
Since $\pi(x;q,a)$ does not change if we shift $a$ by an integer multiple of $q$, we can assume 
$0 \leq a \leq q-1$. Moreover,  if $\gcd(a,q) > 1$, then $\pi(x;q,a) \leq 1$ and the lemma holds. 

So, we can assume $\gcd(a,q)=1$.

Let $S$  be the set of positive integers not exceeding $x$ which are congruent to $a$ modulo $q$. 
Note that $\pi(x;q,a)$ is the number of elements of $S$ which are prime.

Denote the size of $S$ by $k$, that is $k:=|S|$. Then, $k = \lfloor \frac{x-a}{q} \rfloor + 1$. 
We consider six cases depending on what is the remainder when $k$ is divided by $6$. 

Case 1. $k=6l$.

Here we  divide the elements of $S$ into $l$ sextuples of consecutive elements,  $(a,  \ldots,  a+5q), \ldots, (a+(6l-6)q,\ldots, a+(6l-1)q)$. 
Since $\gcd(a,q)=1$ in each sextuple we have a complete set of residues modulo $6$. Therefore, in each sextuple we have exactly two elements which are relatively prime to $6$. 
So, in each sextuple (except possible the first one) there are at most two primes. The first sextuple may contain three primes (if $a=2$ or $a=3$) but not four. 
So, if $k=6l$,  $\pi(x;q,a) \leq 2l+1$. 

Case 2. $k=6l+1$

Here we put all elements of $S$, except the first one into $l$ sextuples of consecutive elements. None of the sextuples contain $2$ or $3$, so 
 if $k=6l+1$, $\pi(x;q,a) \leq 2l+1$. 
 
Case 3. $k=6l+2$

Here we put all elements of $S$, except the first two into $l$ sextuples of consecutive elements.  Again,  none of the sextuples contain more than two primes, so 
 if $k=6l+2$, $\pi(x;q,a) \leq 2l+2$.
 
Case 4. $k=6l+3$

Here we put all elements of $S$, except the first three into $l$ sextuples of consecutive elements. Again,  none of the sextuples contain more than two primes. Moreover, the first three elements of $S$ which are 
$a, a+q,a+2q$  are not all prime since either $a+q$ or $a+2q$ is even and both are greater than $2$. Therefore,
 if $k=6l+3$, $\pi(x;q,a) \leq 2l+2$.
  
Case 5. $k=6l+4$

Here we put all elements of $S$, except the first four into $l$ sextuples of consecutive elements.Again,  none of the sextuples contain more than two primes. Moreover, the first four elements of $S$ cannot be all prime, since as in the previous case either $a+q$ or $a+2q$  is composite. Therefore,
 if $k=6l+4$, $\pi(x;q,a) \leq 2l+3$.
 
 Case 6.  $k=6l+5$

Here we put all elements of $S$, except the first five into $l$ sextuples of consecutive elements.  Here we claim that among the first five 
elements of $S$, $a,a+q,a+2q,a+3q,a+4q$ there are at most three primes. Indeed, if $a$ is even, then $a+2q$ and $a+4q$ are composite (both even and greater than $2$). If $a$ is odd, then 
$a+q$ and $a+3q$ are composite (even and greater than $2$). So,  if $k=6l+4$, $\pi(x;q,a) \leq 2l+3$.
 
 Thus, in all six cases $$\pi(x;q,a) \leq \frac{k}{3} + \frac{5}{3}.$$

Recalling that $k = \lfloor \frac{x-a}{q} \rfloor + 1$ completes the proof of the lemma. 
\end{proof}

Using the above lemma,  we obtain

\begin{equation} \label{c31} 
  Q_3 \leq \sideset{}{'} \sum_{ n^{1/8}  < q \leq n^{1/4}} \pi(\sqrt{n};q^2,n) \leq     \sideset{}{'} \sum_{ n^{1/8}  < q \leq n^{1/4}} \left (  \frac{\sqrt{n}}{3q^2} + 2\right )
\end{equation}

Define
$$c_3(A,B) := \sideset{}{'} \sum_{A < q \leq B} \frac{1}{q^2}.$$

Then, 
\begin{equation} \label{Q31}
Q_3 \leq \frac{c_3\left (n^{1/8},n^{1/4}\right )\sqrt{n}}{3} + 2\pi(n^{1/4}) - 2\pi(n^{1/8})
\end{equation}

 To estimate $c_3(A,B)$ we use the following lemma.

 \begin{lemma} \label{q^2} 
  Suppose $A > 2$   and 
  \begin{equation} \label{uplow}
  \frac{c_4t}{\log t} \leq \pi(t) \leq \frac{c_5t}{\log t}  
  \end{equation}
  for some constants $c_5 > c_4 > 0$ and all $t \geq A$. 
  Then, 
  \begin{equation} \label{q3}
 \sideset{}{'}   \sum_{ q  < A}  \frac{1}{q^2} < \frac{2c_5 - c_4}{A\log A}.
  \end{equation}
\end{lemma}
 
 \begin{proof}
 Using a Riemann-Stieltjes integral we have 
 $$  S(A):= \sideset{}{'} \sum_{   A < q}  \frac{1}{q^2} = \int_A^\infty \frac{d(\pi (t))}{t^2}.$$
 Next, integrating by parts we get 
 $$S(A) = \left.\frac{\pi(t)}{t^2}\right \vert_A^\infty + 2\int_A^\infty \frac{\pi (t)}{t^3}dt.$$
 Using that ${\displaystyle  \pi(t) \leq \frac{c_5t}{\log t}}$ for $t \geq A$ we obtain, 
 $$S(A) \leq  -\frac{\pi(A)}{A^2} + 2c_5\int_A^\infty \frac{1}{t^2 \log t}dt.$$
 Since, $\log t$ is increasing function for $t>2$, $$\int_A^\infty \frac{1}{t^2 \log t}dt < \frac{1}{\log A} \int_A^\infty \frac{1}{t^2}dt = \frac{1}{A \log A}.$$
Finally, using that   ${\displaystyle \frac{c_4A}{\log A} \leq \pi(A)}$ we obtain the lemma.
 \end{proof}
  
 Next, we use the bounds of Rosser and Schoenfeld \cite{RS} for $\pi(t)$. 

\begin{theorem} (J.~Rosser and L.~Schoenfeld) \label{RS lemma}
We have 
\begin{equation} \label{RS1}
\pi (x) >\frac{x}{\log x}\left (  1 + \frac{1}{2\log x} \right )  \quad \mbox{for } x \geq 59,
\end{equation}
and
 \begin{equation} \label{RS2}
\pi (x) <\frac{x}{\log x}\left (  1 + \frac{3}{2\log x} \right )  \quad \mbox{for } x >1.
\end{equation}
\end{theorem}

Using Theorem \ref{RS lemma} we get that if $A \geq 59$, we can take $c_5 = 1 + \frac{3}{2\log A}$ and 
 $c_4 = 1 + \frac{1}{2\log A}$. Therefore, for $A \geq 59$ and $B > A$, 
 $$c_3(A,B)  < \frac{1}{A \log A} + \frac{5}{2\log ^2 A}.$$
 
 Therefore, 
 \begin{equation} \label{c3est}
 c_3\left (n^{1/8},n^{1/4} \right ) < \frac{ 8\log n + 160}{n^{1/8}\log ^2 n},
 \end{equation} 
 for $n > 59^8$.
 
 We obtain
 \begin{equation} \label{Q32}
 Q_3 \leq  \frac{\sqrt{n}(8\log n + 160)}{3n^{1/8}\log ^2 n}  + 2\pi(n^{1/4}) - 2\pi(n^{1/8}) 
 \end{equation}

{\bf Case 4.} $n^{1/4} <  q  \leq     \frac{\sqrt{n}}{c\log n}$, where $c$ is a fixed constant in the interval $(0.5,5]$.

 We already noted above that $f_1(n) =\frac{n^{1/4}}{\log n}$   is increasing for $n > e^4$ and $f(10^{9})=8.55\ldots$.
 Thus,  $ \frac{\sqrt{n}}{10\log n}> n^{1/4}$ for $n \geq 10^{9}$.
 
Let $Q_4$ be the number of primes $p \leq  \sqrt{n}$  such that congruence \eqref{cong} holds for some $q \in (n^{1/4},   \frac{\sqrt{n}}{c\log n}]$. 
Thus, 

\begin{equation} \label{Q41}
Q_4  \leq  \sideset{}{'} \sum_{n^{1/4} < q \leq   \frac{\sqrt{n}}{c\log n}}  \  \pi(\sqrt{n};q^2,n).
\end{equation}

Note that if $q \geq n^{1/4}$, then $\pi(\sqrt{n};q^2,n) \leq 1$  since there will be at most one positive integer congruent to $n$ modulo $q^2$ not exceeding $\sqrt{n}$. Thus, 

\begin{equation} \label{1/8 1/4}
Q_4  \leq \pi \left (   \frac{\sqrt{n}}{c\log n} \right ) - \pi \left ( n^{1/4} \right ).
\end{equation}

Denote $m  =\frac{\sqrt{n}}{c\log n}$. We proved above that $m > n^{1/4}$ for $n>10^{12}$, so $\log m >\log n /4$ for such $n$.  By Theorem
\ref{RS lemma} we have 
$$\pi (m) < \frac{m}{\log m} < \frac{4\sqrt{n}}{c\log^2 n}\left ( 1 + \frac{6}{\log n} \right ),$$
and
$$\pi \left ( n^{1/4}  \right )  < \frac{4n^{1/4}}{\log n}\left ( 1 + \frac{6}{\log n} \right ).$$

Combining  the above two inequalities with equations \eqref{Q32} and \eqref{1/8 1/4} we obtain
\begin{equation} \label{c33}
Q_3+Q_4 < \frac{\sqrt{n}(8\log n + 160)}{3n^{1/8}\log ^2 n} +  \frac{4\sqrt{n}}{c\log^2 n}\left ( 1 + \frac{6}{\log n} \right ) +  \frac{4n^{1/4}}{\log n}\left ( 1 + \frac{6}{\log n} \right ).
\end{equation}
 
{\bf Case 5.}  $\sqrt{n}/(c \log n) <  q < \sqrt{n}$

Let $Q_5$ be the number of primes $p \leq  \sqrt{n}$  such that congruence \eqref{cong} holds for some $q  \in (\frac{\sqrt{n}}{c\log n}, \sqrt{n}]$. 
Thus, 

\begin{equation} \label{Q51}
Q_5  \leq  \sideset{}{'} \sum_ {\frac{\sqrt{n}}{c\log n}  < q < \sqrt{n}} \  \pi(\sqrt{n};q^2,n).
\end{equation}

To estimate $Q_5$  we use the following.
  
  
  Suppose $q_1$ and $q_2$ are distinct primes in $(\sqrt{n}/(c \log n), \sqrt{n})$ such that there exist distinct primes $p_1, p_2$ both not exceeding $\sqrt{n}$ such that $q_1^2 \mid n-p_1$ and $q_2^2 \mid n-p_2$. Thus, there exist positive integers 
  $k_1$ and $k_2$ such that $n-p_1 = k_1q_1^2$ and  $n-p_2=k_2q_2^2$. 
  
  We claim that $k_1 \neq k_2$. Assume the opposite, that $k_1=k_2$
  
  Then, $(n-p_1) - (n-p_2) = k_1(q_1^2 - q_2^2)$, that is
  
 \begin{equation} \label{k1k2} 
     p_2 - p_1 = k_1(q_1+q_2)(q_1-q_2).  
  \end{equation} 
  
  Note that $k_1q_1 = k_1q_1^2/q_1 = (n-p_1)/q_1 > (n-\sqrt{n})/ \sqrt{n} = \sqrt{n}-1$.   
  
  Therefore, $k_1(q_1+q_2)|q_1-q_2| > k_1q_1 > \sqrt{n}-1  > |p_1 - p_2|$,  contradicting
  \eqref{k1k2}. So, $k_1 \neq k_2$. 
  
  We conclude that the contribution from $q \in (\sqrt{n}/(c \log n), \sqrt{n})$ does not exceed the number of distinct integers $k$ such that $n- p = kq^2$. Note that for each such $k$ we have $k = (n-p)/q^2 < (c \log n)^2$, since 
  $q > \sqrt{n}/(c \log n)$. Thus, 
  
  \begin{equation} \label{c4}
Q_5  \leq  (c \log n)^2.
  \end{equation}

Now, combining \eqref{Q13}, \eqref{r21}, \eqref{c33}, and \eqref{c4} we get 
 $$
   Q_1+Q_2+Q_3+Q_4+Q_5 <   \frac{46\sqrt{n}}{35\log n} +  \frac{1569\sqrt{n}}{350\log^2 n}
   +    \frac{4C_1(11)\sqrt{n}}{\log n}+  \frac{32C_2(11)\sqrt{n}}{\log^2 n}+$$
 $$+   \frac{\sqrt{n}(8\log n + 160)}{3n^{1/8}\log ^2 n} +  \frac{4\sqrt{n}}{c\log^2 n}\left ( 1 + \frac{6}{\log n} \right ) +  \frac{4n^{1/4}}{\log n}\left ( 1 + \frac{6}{\log n} \right )+ (c \log n)^2,$$
 for $n > 59^8$. 
 
 Since to prove the theorem it is sufficient to show $Q_1+Q_2+Q_3+Q_4+Q_5  < \pi \left ( \sqrt{n} \right )$ and by Theorem \ref{RS lemma} for $n > 59^2$, 
 $$\pi \left ( \sqrt{n} \right ) > \frac{2\sqrt{n}}{\log n}\left (1 + \frac{1}{\log n}\right ),$$ it is sufficient to establish 
 
 \begin{equation} \label{c5}
 \frac{2\sqrt{n}}{\log n}\left ( 1 + \frac{1}{\log n} \right )  > f_2(n), 
 \end{equation}
 where 
 $$f_2(n) :=   \frac{46\sqrt{n}}{35\log n} +  \frac{1569\sqrt{n}}{350\log^2 n}
   +   \frac{4C_1(11)\sqrt{n}}{\log n}+ \frac{32C_2(11)\sqrt{n}}{\log^2 n}+$$
 $$   +\frac{\sqrt{n}(8\log n + 160)}{3n^{1/8}\log ^2 n} +  \frac{4\sqrt{n}}{c\log^2 n}\left ( 1 + \frac{6}{\log n} \right ) +  \frac{4n^{1/4}}{\log n}\left ( 1 + \frac{6}{\log n} \right )+ (c \log n)^2.$$

Dividing the   inequality \eqref{c5} by $\frac{\sqrt{n}}{\log n}$ and simplifying we get that the inequality is equivalent to 
\begin{equation} \label{c6}
\frac{24}{35} > \frac{869}{350\log n} + 4C_1(11) + \frac{32C_2(11) }{\log n}+\frac{ 8\log n + 160}{3n^{1/8}\log  n} + 
\end{equation}
$$+\frac{4 }{c\log  n}\left ( 1 + \frac{6}{\log n} \right )+ 4n^{-1/4} \left ( 1 + \frac{6}{\log n} \right )+ \frac{c^2 \log ^3 n}{\sqrt{n}}.$$

Note the functions $\frac{1}{\log n}$, $\frac{1}{\log^2 n}$, $\frac{1}{n^{1/8}}$, $\frac{1}{n^{1/8}}$, $\frac{1}{n^{1/8}\log n}$, and $\frac{1}{n^{1/4}\log n}$ are all decreasing for $n > 1$. Moreover, the function 
$\frac{\log^3 n}{\sqrt{n}}$ is decreasing for $n > e^6$. Therefore, if inequality \eqref{c6} holds for some $n_0 > 59^8$, then it holds for all $n \geq n_0$. Since the inequalities \eqref{c5} and \eqref{c6} are equivalent, 
then if  inequality \eqref{c5} holds for some $n_0 > 59^8$, then it holds for all $n \geq n_0$.

Estimating the tail of the series defining $C_1(11)$ and $C_2(11)$ in the same way we estimated the series   $\sideset{}{'} \sum_{A  < q < B} \frac{1}{q^2}$ we get the estimates $C_1(11)<0.033$ and $C_2(11)<0.1$. 

Taking $c=4$ and $n=59^8 + 1$, we get that the left-hand-side of inequality \eqref{c5} exceeds the right-hand-side by more than $95945 > 0$. Thus, the theorem holds for $n > 59^8$. 

From now on, we assume $n \leq 59^8$.  For such $n$ we cannot use Theorem \ref{RS lemma} to get lower bound for $\pi \left (n^{1/8} \right )$. However, we can estimate $c_3(A,B)$ as follows. 
We have
$$c_3(A,B) < c_3(A) =  \sideset{}{'} \sum_q \frac{1}{q^2} -  \sideset{}{'} \sum_{q \leq A} \frac{1}{q^2}.$$

Denote $$g(A) :=  \sideset{}{'} \sum_{q \leq A} \frac{1}{q^2}.$$
The function $g$ is piecewise continuous and nondecreasing. 

Also, it has been known since Euler \cite{Euler} p.480 (or possibly earlier)  that 
\begin{equation} \label{1/q^2}
 \sideset{}{'} \sum_q \frac{1}{q^2} := c_0 =0.452247420041065\ldots
\end{equation} 

Replacing the estimate for $c_3 \left ( n^{1/8},n^{1/4} \right )$ in   equation \eqref{c33} we get 
\begin{equation} \label{c331}
Q_3+Q_4 < \frac{\sqrt{n}\left (c_0  - g\left ( n^{1/8}\right )\right )}{3} +  \frac{4\sqrt{n}}{c\log^2 n}\left ( 1 + \frac{6}{\log n} \right ) +  \frac{4n^{1/4}}{\log n}\left ( 1 + \frac{6}{\log n} \right ).
\end{equation}

To prove the theorem now we need to establish 
\begin{equation} \label{c51}
 \frac{2\sqrt{n}}{\log n}\left ( 1 + \frac{1}{\log n} \right )  > f_3(n), 
 \end{equation}
 where 
 $$f_3(n) :=   \frac{46\sqrt{n}}{35\log n} +  \frac{1569\sqrt{n}}{350\log^2 n}
   +   \frac{4C_1(11)\sqrt{n}}{\log n}+ \frac{32C_2(11)\sqrt{n}}{\log^2 n}+$$
 $$   +\frac{\sqrt{n}\left (c_0  - g\left ( n^{1/8}\right )\right )}{3} +  \frac{4\sqrt{n}}{c\log^2 n}\left ( 1 + \frac{6}{\log n} \right ) +  \frac{4n^{1/4}}{\log n}\left ( 1 + \frac{6}{\log n} \right )+ (c \log n)^2.$$

Unfortunately, if we divide $f_3(n)$ by $\frac{\sqrt{n}}{\log n}$ it is not clear that the resulting function is decreasing due to the term $  \frac{\sqrt{n}\left (c_0  - g\left ( n^{1/8}\right )\right )}{3}$, so we proceed as follows. 

Let $p_k$ be the $k$th prime number. Next we will be considering the cases when $n \in (p_k^8, p_{k+1}^8]$ where $31 \leq p_k \leq 53$ (when $p_k=31$ we will consider $n \in [10^{12}, 37^8]$). For such $n$
we have $g(n^{1/8}) \geq g(p_k)$. Moreover, since $\log n$ is increasing function, $\frac{8\log p_{k+1}}{\log n} \geq 1$ for $n \in (p_k^8, p_{k+1}^8]$. Also, since now we have upper bound for $n$, when estimating $Q_2$ we be using equation \eqref{r2}, that is we will use $c_1(11,p_{k+1})$ instead of $C_1(11)$, and $c_2(11,p_{k+1})$ instead of $C_2(11)$.

We replace $f_3(n)$ by the (potentially) slightly larger function 
 $$f_4(n,k) := \frac{46\sqrt{n}}{35\log n} +  \frac{1569\sqrt{n}}{350\log^2 n}
   +   \frac{4c_1(11, p_{k+1})\sqrt{n}}{\log n}+ \frac{32c_2(11, p_{k+1})\sqrt{n}}{\log^2 n}+$$
 $$   +\frac{\sqrt{n}\left (c_0  - g\left ( p_k \right )\right )(8\log p_{k+1})}{3\log n} +  \frac{4\sqrt{n}}{c\log^2 n}\left ( 1 + \frac{6}{\log n} \right ) +  \frac{4n^{1/4}}{\log n}\left ( 1 + \frac{6}{\log n} \right )+ (c \log n)^2.$$
The function $f_4(n,k)$ may be   larger than $f_3(n)$ but has the property that it is decreasing on the interval $(p_k^8, p_{k+1}^8]$ after division by $\frac{\sqrt{n}}{\log n}$. 

Thus, if
\begin{equation} \label{r7}
 \frac{2\sqrt{n_0}}{\log n_0}\left ( 1 + \frac{1}{\log n_0} \right )  > f_4(n_0,k)
 \end{equation}
for some $p_{k+1}^8 \geq n_0 \geq \max(p_k^8, 10^{12})$, then the theorem holds for all $n \in [n_0, p_{k+1}^8]$. 

In the table below we record the results of our computations in several intervals. The symbol $\Delta $ will denote the difference between the left-hand-side of inequality \eqref{r7} and its right-hand-side evaluated at the left end of each interval  indicated in each row.

$$\begin{array}{|c|c|c|c|c|c|c|} \hline
k & \mbox { interval }& c_1(11,p_{k+1}) & c_2(11,p_{k+1}) & c & c_0 - g(p_k) & \Delta\\ \hline
16 & [53^8, 59^8] & 0.02941652 & 0.08277361 & 4 & 0.00352137 & 74613.3\\ \hline
15 & [47^8,53^8] & 0.02912429 & 0.08158205 & 3.3 & 0.00387736 & 46560.4\\  \hline
14 & [43^8,47^8] & 0.02876145 & 0.08014145 & 2.9 & 0.004330053 & 32612.4 \\ \hline
13 & [41^8,43^8] & 0.02829891 &  0.07836062 & 2.7 & 0.00487089 & 26933.3 \\ \hline
12& [37^8,41^8]& 0.02774520  & 0.07627800 & 2.3 & 0.00546577 & 17911.1 \\ \hline
11&[10^{12},37^8]& 0.02713545 & 0.07401364& 2 & 0.00619623 &9029.5 \\ \hline
\end{array}
$$

The left end of the interval in the last row in the  table above is $10^{12}$ rather than $31^8$ since $31^8 < 10^{12}$ and we 
cannot use Theorem \ref{BM} to estimate $\pi(\sqrt{n}; q^2,n)$ when $\sqrt{n} < 10^6$. 

When $n \leq 10^{12}$ we use another result of Bennett et. al., Corollary 1.6 of \cite{BenMar}.

\begin{theorem} \label{BenMar1}
Let $1 \leq q \leq 1200$ be an integer, and $a$ be an integer coprime to $q$. For all $x \geq 50q^2$ we have 
$$\frac{x}{\varphi(q)\log x} <  \pi(x;q,a)  < \frac{x}{\varphi(q)\log x} \left ( 1 + \frac{5}{2\log x} \right ).$$
\end{theorem}

We apply the above theorem with $n \geq 10^9$, and $q \in \{2^2,3^2,5^2,7^2\}$. Since $\sqrt{10^9} > 31622 > 50\cdot 7^2 =2450$, Theorem \ref{BenMar1} applies when $\gcd(n,210)=1$  and argue similarly to Case 1 when $\gcd(n,210)>1$. 
Substituting into equation \eqref{Q11} we obtain

\begin{equation} \label{newQ11}
Q_1 < \frac{46\sqrt{n}}{35\log n} +  \frac{32\sqrt{n}}{7\log^2 n}.
\end{equation}

We get a somewhat worse estimate. The constant in the second term is $32/7 = 4.57\ldots$ rather than $1569/350=4.48\ldots$.

So, for $n \geq 10^9$ now we need to show 

\begin{equation} \label{f5}
 \frac{2\sqrt{n}}{\log n}\left ( 1 + \frac{1}{\log n} \right )  > f_5(n,k),
\end{equation}
where
$$f_5(n,k) := \frac{46\sqrt{n}}{35\log n} +  \frac{32\sqrt{n}}{7\log^2 n}
   +   \frac{4c_1(11, p_{k+1})\sqrt{n}}{\log n}+ \frac{32c_2(11, p_{k+1})\sqrt{n}}{\log^2 n}+$$
 $$   +\frac{\sqrt{n}\left (c_0  - g\left ( p_k \right )\right )(8\log p_{k+1})}{3\log n} +  \frac{4\sqrt{n}}{c\log^2 n}\left ( 1 + \frac{6}{\log n} \right ) +  \frac{4n^{1/4}}{\log n}\left ( 1 + \frac{6}{\log n} \right )+ (c \log n)^2.$$
 
 We proceed exactly as before, the only difference between $f_4(n,k)$ and $f_5(n,k)$ is that the constant $1569/350$ in the second term of $f_4(n,k)$ is replaced by $32/7$ to obtain $f_5(n,k)$. 
 
 The new table with $f_5$ rather than $f_4$ is: 
 
 $$\begin{array}{|c|c|c|c|c|c|c|} \hline
k & \mbox { interval }& c_1(11,p_{k+1}) & c_2(11,p_{k+1}) & c & c_0 - g(p_k) & \Delta\\ \hline
11 & [31^8, 10^{12}] & 0.02713545 & 0.07401364 & 2 & 0.00619623 & 5606.9\\ \hline
10 & [29^8,31^8] & 0.02638469 & 0.0713027 & 2 & 0.00723681 & 3669.5\\  \hline
9 & [23^8,29^8] & 0.02530942 & 0.06761027 & 1.6 & 0.00842587 & 910.8 \\ \hline
8 & [19^8,23^8] & 0.02407789 &  0.0634633 & 1.3 & 0.01031623 & 179.8 \\ \hline
7& [17^8,19^8]& 0.02210161  & 0.05726672 & 1.1 & 0.0130863 & 62.1 \\ \hline
6&[13^8,17^8]& 0.01917763 & 0.04865725& .6 & 0.0165465 &35.3 \\ \hline
\end{array}
$$

Thus, the theorem is holds for $n \geq 13^8 = 815730721$. Since we established the theorem by direct computation for $n \leq 10^9$, this completes the proof.

\end{proof}

\end{document}